\documentclass{article}
\usepackage[T1]{fontenc}
\usepackage[latin1]{inputenc}
\usepackage{setspace}

\makeatletter

\providecommand{\LyX}{L\kern-.1667em\lower.25em\hbox{Y}\kern-.125emX\@}

\makeatother

\begin{document}

\( \qquad  \)\( \qquad  \)A Computational Algorithm for \( \pi  \)( N )\\
\\
\\
 \\
\\
\( \qquad \qquad  \)\( \qquad  \)Abhijit Sen\\
\( \qquad  \)Theory Group,Saha Institute of Nuclear Physics\\
\( \qquad  \)1/AF,BidhanNagar,Calcutta-700064\\
\( \quad  \)\( \quad  \) West Bengal , INDIA\\
\\
\\
\( \qquad \qquad \qquad  \)Satyabrata Adhikari\\
\( \qquad  \)Calcutta Mathematical Society, A.E 374\\
\( \qquad  \)BidhanNagar , Sector-I, Calcutta-700064\\
\( \qquad  \)West Bengal , INDIA\\
\\
\\
\\
\\
\\
\\
\\
\( \qquad \qquad \qquad \qquad \qquad \qquad  \)Abstract\\
\\
\\
\\
\( \qquad  \)An algorithm for computing \( \pi  \)(N) is presented.It is shown
that using a symmetry of odd composites we can easily compute \( \pi  \)(N).This
method relies on the fact that counting the number of odd composites not exceeding
N suffices to calculate \( \pi  \)(N).
AMS classification number(1980):10-01\\
\\
\\
\\
\\
Introduction\\
\\
 \\
\( \qquad  \)The problem of finding a compact mathematical formula for \( \pi  \)(N)
is well known {[}1{]}.Though no compact result has been reported till date,
a number of accurate/fairly accurate algorithms exist {[}2,3,7,8{]} among which
the prime number theorem {[}4{]} and the function Li(x) {[}3{]} require special
mention.Wilson's theorem {[}5{]} seems to be unique,however,like all reportings
so far it has a low practical utility. Recently a result has been reported using
Smarandache Functions{[}6{]} which too however has a low practical utility.\\
\\
\( \qquad  \)In the present paper we develope a computational algorithm based
on a certain symmetry of odd composites.This is the main essence of the present
piece of work.\\
\( \qquad  \)\\
\\
 The Algorithm\\
 \\
A. General Formulation\\
\\
\( \qquad  \)Among all natural numbers we know that \hfill{}\\

i) All even numbers except 2 are composite.

ii)All odd composites can be factorised into odd primes only.

iii)Unity is neither prime nor composite.\\
\\
Therefore we have

\( \pi  \)(N) = N - n\( _{e} \)(c) - n\( _{o} \)(c) - 1 \( \qquad  \)...........
(1)\\
where

n\( _{e} \)(c) and n\( _{o} \)(c) give respectively the number of even and
odd composite numbers not exceeding N.Without the loss of generality we assume
N\( \geq  \)2.\\
\\
Now we have\\
\( \qquad n_{e} \)(c) = Floor (\( \frac{N}{2} \)- 1 ) \( \qquad  \).............(1a)\\
\\
To obtain n\( _{o} \)(c) let us refer to the multiplication table of odd natural
numbers \( \geq  \)3 as shown in Table-II. It shows a very small section at
the top left hand corner of the infinite multipliction table.\\
 \\
It can be shown that all odd composites occur at least once within Table -II
whose general term is given by \\
t\( _{n,r} \) = ( 2n + 1 )\( ^{2} \)+ 2 ( r - 1 ) ( 2 n + 1 ) = ( 2 n + r
)\( ^{2} \)- ( r - 1 )\( ^{2} \)....................(2)\\
 where both n and r are natural.\\
We further see that for each 'n' value we have an A.P originating at ( 2n +
1)\( ^{2} \)\\
with a common difference 2 ( 2 n + 1 ). The terms of the A.P for each 'n' are
\\
assumed to constitute a set represented as \\
\{ t\( _{n,r} \) \} = \{ x : x = (2 n + r )\( ^{2} \)- ( r - 1 )\( ^{2} \)\}
\\
where r varies from 1 to 

r|\( _{max} \) = Floor (\( \frac{1}{2(2n+1)} \)( N - ( 2n+1)\( ^{2} \)) +
1 ) .................. (3a)\\
for a given n i.e for a given set \{ t\( _{n,r} \) \}. \\
\\
Thus,\\
 n\( _{o} \)( c ) = \( \sum  \) Floor (\( \frac{N-(2n+1)^{2}}{2(2n+1)} \)
+ 1 ) - \( \lambda _{c} \) ..................(4a)\\
where the sum is over n from 1 to 

n|\( _{max} \) = Floor (\( \frac{\sqrt{N}-1}{2} \)) ......... (3b). 

Here \( \lambda _{c} \) is a correction term which we would discuss shortly.\\
\\
\\
 B.Corrections\\
 \\
\( \qquad  \) \\
B1. The Summation\\
\\
\\
\( \quad  \) For a composite (2n\( _{r} \)+1),all the terms arising from the
summation in (4a) corresponding to n=n\( _{r} \) would be a duplication of
some previously occuring term for some n = n\( _{s} \) < n\( _{r} \) .Therefore
we shall restrict ourselves to only those values of n\( _{r} \) for which (2n\( _{r} \)
+ 1) is prime.\\
\( \qquad  \)Following the discussions as above we can recast the summation
part of (4a) as\\
\\
n\( _{o} \)(c) = \( \sum _{i}^{R} \) (Floor\( \frac{N-(2n_{i}+1)^{2}}{2*(2n_{i}+1)} \)+1)-\( \lambda _{c} \)...(4b)\\
where we restrict the sum for a prime (2n\( _{i} \)+1). Here the suffix R denotes
the sum over the restricted sets of the values of n's.The primality of (2n\( _{i} \)+1)
can be ascertained using the symmetry we mentioned.We do not accept those values
of n\( _{k} \) for which 

\( \frac{1}{2} \)( t\( _{n,r} \) -1 ) = n\( _{k} \) \( \leq  \)n\( _{max} \)
for r \( \leq  \) r\( _{max} \) as given by (3a) and (3b).\\
As will be shown in a latter article {[} 9 {]} all composite (2n\( _{i} \)
+ 1 ) can be expressed as a difference of two squares following eqn. (2) .

Alternatively the primality could be determined using Wilson's theorem{[}5{]}.\\

B2. Evaluation of \( \lambda _{c} \)\\
\\
\( \qquad  \)The need of \( \lambda _{c} \) arises from the fact that for
different n values the sets \{ t\( _{n,r} \)\} coincide at certain points.\\
\\
 Now once the summation in (4a) be restricted we shall also require to restrict
the evaluation of \( \lambda _{c} \).\\
 \\
\( \qquad  \)Restricting ourselves to those values of n\( _{r} \) for which
(2n\( _{r} \)+1) is prime we can evaluate \( \lambda _{c} \) using the set
theoretic formula\\
\\
n(A\( _{1} \)\( \bigcup  \)A\( _{2} \)\( \bigcup  \)....\( \bigcup  \)A\( _{k} \))
= \( \sum  \) n(A\( _{i} \))-\( \sum  \)\( \sum  \)n(A\( _{i} \)\( \bigcap  \)A\( _{j} \))+\( \sum  \)
n(A\( _{i} \)\( \bigcap  \)A\( _{j} \)\( \bigcap  \)A\( _{k} \)) - .......
(5)\\
where n(X) denotes the cardinality of set X.\\
Here summations over all indices can range from 1 to n\( _{max} \) with the
restriction\\
i < j < k ... and so on.\\
Let us define the following notations:\\
For intersection of two A.P's let \\
t\( _{n_{1},n_{2}}^{(1)} \) be the first common term of \{ t\( _{n_{1},r} \)\}
and \{ t\( _{n_{2},r} \) \} while\\
( c.d )\( _{n_{1},n_{2}} \) be the common difference of the common terms \\
\( \quad  \)\( \quad  \) For intersection of more than two A.P's the notation
is identical.\\
\( \qquad  \)It can be shown that\\
\\
a ) For two A.P's\\
\\
t\( _{n_{1},n_{2}}^{(1)} \)= (2n\( _{1} \)+1 ){*} ( 2n\( _{2} \) + 1 ) {*}
\{ 3 + 2 {*} Ceil {[} \( Z \) {]}{*} H {[} Z {]} \} ... (6a)\\
\\
where\\
Z = \( \frac{1}{2} \)(\( \frac{2n_{2}+1}{2n_{1}+1} \) - 3 ) \\
and\\
( c . d )\( _{n_{1},n_{2}} \)= 2 {*} ( 2n\( _{1} \)+ 1 ) {*} ( 2 n\( _{2} \)
+ 1 ) ... (6b)\\
Here H {[} Z {]} denotes the Heaviside Step Function defined by\\
\\
H{[}z{]}= 1 , z > 0

=\( \frac{1}{2} \), z=0

= 0 , z < 0\\

b ) For more than two A.P's\\
\\
t\( _{n_{1},n_{2}...n_{k}}^{(1)} \)= \( \prod  \)\( _{i=1}^{k} \) ( 2n\( _{i} \)+
1 ) {*} ( 1 + 2 {*} Ceil ( Q) {*} H {[} Q {]} ) ...(7a)\\
\\
where\\
\\
Q = \( \frac{1}{2} \)(\( \frac{(2n_{k}+1)^{2}}{\prod _{i=1}^{k}(2n_{i}+1)} \)-
1 )\\
and\\
( c.d )\( _{n_{1}},n_{2},...,n_{k} \) = 2 {*} \( \prod ^{k}_{i=1} \)( 2n\( _{i} \)
+ 1 ) ...(7b) \\
\\
Using (5),(6a),(6b),(7a) and (7b) we get \\
\\
\( \lambda _{c} \)=\( \sum ^{n_{max}}_{r=2} \)(-1)\( ^{r} \)\( \sum  \)\( ^{R} \)
\( \sum  \)\( ^{R} \)...\( \sum ^{R} \)(Floor(B)+1){*}H{[}B+\( \varepsilon  \){]}
...(8)\\
\\
where\\
i)The summations run over all restricted sets of n values,one for each set with
\\
n\( _{1} \) < n\( _{2} \) < n\( _{3} \) ..... < n\( _{k} \)\\
ii)\( \varepsilon  \) is a positive number,however small.It will be neglected
at the last stages.\\
iii)We have \\
\\
 B = \( \frac{N-t_{n_{1},n_{2},...,n_{r}}^{(1)}}{(c.d)_{n_{1},n_{2},...,n_{r}}} \) 

The Heaviside function restricts the numerator to non-negative values.

Using (1),(1a),(4b),(5), and (8) the result follows.\\
Though equation (8) seems to have many steps in practice it would be quite limited
as once the result given by (6a) and (7a) exceeds N,higher n\( _{r} \) terms
would be in general very very limited in number.\\
\\
The Results\\
\\
\( \qquad  \)The results as obtained from the above mentioned algorithm are
given in Table I.It is seen that the value of \( \pi  \)(N) is exactly reproduced.This
establishes the correctness of the above mentioned algorithm.\\
\( \quad \quad  \)To illustrate the method let us evaluate \( \pi  \)(100).We
have\\
from (3b) : n\( _{max} \)=4\\
Now,to calculate \( \pi (100) \) we need to know all primes till \( \sqrt{100} \)=10.\\
Corresponding to n=4 we get (2n + 1)=9 < 10 which is a composite.\\
So we do not consider n=4.\\
Proceeding:\\
from ( 1a ): n\( _{e} \)(c)=Floor(\( \frac{100}{2}-1) \) = 49\\
from (4b):n\( _{o} \)(c)= 16+8+4=28\\
from ( 8 ):\( \lambda _{c} \)=2+1=3\\
Therefore from (1) we have \( \pi  \)(100) = 100 - 49 - 28 + 3 - 1 = 25\\

Conclusion\\
\\
\( \qquad  \)The formula for \( \pi  \)(N) as established seems complicated
at first sight.However the symmetry of natural numbers used, the necessity of
introduction of \( \lambda  \)\( _{c} \) and the simple set-theoretic prescription
for its evaluation makes the present method quite tractable. However, like the
well known {}``Sieve algorithm{}'' {[}8{]} this method too requires the knowledge
of all primes till \( \sqrt{N} \) to evaluate \( \pi (N) \) .The number of
steps is much reduced as we use the symmetry of natural numbers throughout,
first to determine primality of (2n\( _{i} \)+1) if the same is not already
known and afterwards for the final evaluation of \( \pi  \)(N) in comparison
to the Sieve Algorithm.Refering back to the symmetry of odd composites that
we had mentioned,the symmetry is more clearly evident if we use the binary base
of representation specially in Table II.This is evident from the fact that the
common difference of the common differences in Table II equals base squared.
It is believed that this could be a suitable method to determine \( \pi  \)(N).\\
\\
Acknowledgement\\
\\
\( \qquad  \)The authors would like to acknowledge the constant encouragement
by Mr.Samit De of SINP. \\
\\
\\
\\
\\
\\
\\
\\
\\
 \\
\\
 \\

TABLE : I\\
\\
\( \qquad \qquad \qquad  \)Results of determination of \( \pi  \)(N)\\

\vspace{0.3cm}
{\centering \begin{tabular}{|c|c|c|}
\hline 
N&
\( \pi  \)(N)&
Estimated value\\
\hline 
\hline 
10&
4&
4\\
\hline 
100&
25&
25\\
\hline 
1000&
168&
168\\
\hline 
\end{tabular}\par}
\vspace{0.3cm}

TABLE : II\\
\\
\( \qquad \qquad \qquad  \)Multiplication Table of odd natural numbers\\

\begin{singlespace}
\vspace{0.3cm}
{\centering \begin{tabular}{|c|c|c|c|c|c|c|}
\hline 
&
3&
5&
7&
9&
11&
Common Difference\\
\hline 
\hline 
3&
9&
15&
21&
27&
33&
6\\
\hline 
5&
&
25&
35&
45&
55&
10\\
\hline 
7&
&
&
49&
63&
77&
14\\
\hline 
Common Difference&
&
&
&
&
&
4 =2\( ^{2} \) \\
\hline 
\end{tabular}\par}
\vspace{0.3cm}
\end{singlespace}

\( \qquad \qquad \qquad \qquad  \)\\
\\
\\
\\
\\
\\
\\
\\
\\
\\
\\
 \\
\\
\\
\\
\\
\\
\\
\\
\\
\\
 \\
\\
References\\
\\
\\
1. A-17, Unsolved problems in Number Theory,Richard K Guy, Springer- Verlag,N.Y
\\
2.CRC Concise Encyclopedia of Mathematics, Eric W Weisstein,1630, Chapman \&
Hall/CRC.\\
3.CRC Concise Encyclopedia of Mathematics, Eric W Weisstein,1097, Chapman \&
Hall/CRC.\\
4. Elementary Number Theory ,David M.Burton,Appendixes,396-403(1995)\\
5.Math.Gaz.48(1964)413-415\\
6.CRC Concise Encyclopedia of Mathematics, Eric W Weisstein,1660 - 1661, Chapman
\& Hall/CRC.\\
7. CRC Concise Encyclopedia of Mathematics, Eric W Weisstein,1427 - 1428, Chapman
\& Hall/CRC.\\
8. Elementary Number Theory ,David M.Burton,60-61(1995)\\
9. http://www.arXiv.org, math.GM/0202269
\end{document}